\documentclass[10pt,a4paper,twoside,fleqn,leqno]{article}

\usepackage[utf8]{inputenc}
\usepackage[T1]{fontenc}
\usepackage[full]{textcomp}
\usepackage{amsmath,amssymb}
\usepackage{mathrsfs}

\usepackage{geometry}

\usepackage{mathtools}

\usepackage{morphism}
\usepackage{widebar}

\DeclarePairedDelimiterX{\abs}[1]\lvert\rvert{\ifblank{#1}{\,\cdot\,}{#1}}

\newcommand\SetSymbol[1][]{\nonscript\:#1\vert\allowbreak\nonscript\:\mathopen{}}                                         
\providecommand\given{} % to make it exist                                                                                
\DeclarePairedDelimiterX\Set[1]\{\}{\renewcommand\given{\SetSymbol[\delimsize]}#1}  
\usepackage[final]{microtype}
\usepackage{csquotes}
\usepackage[english]{babel}

\makeatletter
\let\mytagform@=\tagform@
\def\tagform@#1{\maketag@@@{\hbox{{\ignorespaces#1\unskip\@@italiccorr\hspace{1mm}}}}\kern1sp}
\makeatother

\usepackage[thmmarks,amsmath]{ntheorem}
\usepackage{thmtools}

\usepackage{chngcntr}
\counterwithin{equation}{section}

\declaretheoremstyle[headformat=swapnumber,headpunct={.\ ---},headfont=\normalfont\scshape,bodyfont=\itshape,%
spaceabove=0pt,spacebelow=0pt,preheadhook={\medskip}]{theorem}
\declaretheorem[style=theorem,sibling=equation]{theorem}

\declaretheorem[style=theorem,sibling=equation]{lemma}
\declaretheorem[style=theorem,sibling=equation]{corollary}

\declaretheoremstyle[headformat=swapnumber,headpunct={.\ ---},headfont=\scshape,bodyfont=\normalfont,%
spaceabove=0pt,spacebelow=0pt,preheadhook={\medskip}]{definition}
\declaretheorem[style=definition,sibling=equation]{definition}

\declaretheorem[style=definition,sibling=equation]{remark}

\declaretheoremstyle[headpunct=.,headfont=\itshape,bodyfont=\normalfont,%
qed=\ensuremath{\square},spaceabove=0pt,spacebelow=0pt]{nonumberproof}
\declaretheorem[numbered=no]{proof}

\usepackage{booktabs}

\usepackage{enumitem}
\newlist{condition}{enumerate}{5}
\newlist{case}{enumerate}{5}
\setlist{nosep,listparindent=\parindent}
\setlist[itemize]{label=\guillemotright}
\setlist[case,1]{label=\arabic*.,ref=\theequation.\arabic*}
\setlist[case,2]{label=\arabic*.,ref=\thecasei.\arabic*}

\usepackage{listings}

\usepackage{tikz-cd}
\usepackage{lastpage}
\usepackage{fancyhdr}
\usepackage[labelfont=sc]{caption}

\usepackage[
	sorting=nyt,
	doi=false,
	url=false,
	isbn=false,
	backend=bibtex,
	safeinputenc,
]{biblatex}
\bibliography{paper.bib}
\defbibheading{bibliography}[\bibname]{\sectionstar{#1}}
\newcommand\quelle[1]{{%
		\null\nobreak\hfil\penalty50
		\hbox{}\nobreak\hfil{#1}%
		\parfillskip=0pt \finalhyphendemerits=0 \par}}
\renewbibmacro*{pageref}{%
	\iflistundef{pageref}
	{}
	{%
		\printtext{\addperiod\addspace\quelle{\mbox{\scriptsize\printtext[parens]{%
				\ifnumgreater{\value{pageref}}{1}
				{\bibstring{backrefpages}\ppspace}
				{\bibstring{backrefpage}\ppspace}%
			\printlist[pageref][-\value{listtotal}]{pageref}}}}}}}
\usepackage{multicol}

\usepackage{graphicx}
\usepackage{aliascnt}
\usepackage{cleveref}
\creflabelformat{equation}{#2#1#3}
\crefname{akaparagraph}{paragraph}{paragraphs}
\Crefname{akaparagraph}{Paragraph}{Paragraphs}
\crefformat{akaparagraph}{\S#2#1#3}
\crefname{conditioni}{condition}{conditions}
\Crefname{conditioni}{Condition}{Conditions}
\crefname{casei}{case}{cases}
\Crefname{casei}{Case}{Cases}
\crefname{caseii}{case}{cases}
\Crefname{caseii}{Case}{Cases}
\crefname{nproof}{proof}{proofs}
\Crefname{nproof}{Proof}{Proofs}

\newaliascnt{akaparagraph}{equation}

\makeatletter
\newcommand{\sectionstar}[1]{\bigskip%
	{\normalfont\large\scshape \noindent#1}\par\noindent\ignorespaces}
\makeatother

\usepackage{titlesec}
\titleformat{\section}
{\normalfont\large\scshape}
{\thesection}{1em}{}
\titlespacing*{\section}
{0pt}{3.5ex plus .1ex minus .2ex}{1.5ex minus .1ex}

\usepackage{ifthen}
\renewcommand{\paragraph}[1]{\medskip\pagebreak[3]\refstepcounter{akaparagraph}%
	{\normalfont\normalsize\scshape \noindent\theakaparagraph%
	\ifthenelse{\equal{#1}{}}%
	{}%
	{\ #1.}%
	\ ---}%
}

\def\tr{\textsc{tr}}
\def\cm{\textsc{cm}}

\def\HH{\mathrm{H}}
\def\ZZ{\mathbb{Z}}
\def\QQ{\mathbb{Q}}
\def\QQl{\mathbb{Q}_{\ell}}
\def\RR{\mathbb{R}}
\def\CC{\mathbb{C}}

\DeclareMathOperator{\Res}{Res}
\DeclareMathOperator{\Hom}{Hom}
\DeclareMathOperator{\End}{End}
\DeclareMathOperator{\Nm}{Nm}
\DeclareMathOperator{\GL}{GL}
\DeclareMathOperator{\SL}{SL}
\DeclareMathOperator{\SO}{SO}
\DeclareMathOperator{\SU}{SU}
\DeclareMathOperator{\U}{U}

\def\GB{\mathrm{G}_{\textnormal{B}}}
\def\HB{\HH_{\mathrm{B}}}
\def\Gl{\mathrm{G}_{\ell}}
\def\Glc{\mathrm{G}_{\ell}^{\circ}}

\def\scrG{\mathscr{G}}
\def\Hl{\HH_{\ell}}

\def\MTC{\textnormal{MTC}}

\def\HdgB{\textnormal{Hdg}_{\textnormal{B}}}
\def\Hdgl{\textnormal{Hdg}_{\ell}}

\DeclareMathOperator{\reldim}{reldim}
\DeclareMathOperator{\Spec}{Spec}
\DeclareMathOperator{\Gal}{Gal}
\DeclareMathOperator{\Lie}{Lie}
\DeclareMathOperator{\trdeg}{trdeg}
\DeclareMathOperator{\charpol}{charpol}
\def\Frob{\textnormal{Frob}}
\def\der{\textnormal{der}}

\def\tra{\textnormal{tra}}
\DeclareMathOperator{\trace}{tr}

\def\HA{M_{A}}
\def\HX{M_{X}}

\def\title{The Mumford--Tate conjecture for the product of an abelian surface and a K3~surface}
\def\author{J.M.~Commelin}
\usepackage{datetime}
\def\date{\dayofweekname{\day}{\month}{\year}, the \ordinaldate{\day} of \monthname, \number\year}

\begin{document}
\noindent\textsc{\Large \title}

\medskip

\noindent\textit{by} \quad \author \hfill \date

\bigskip

\section{Introduction}

The main result of this paper is the following theorem.
In the next paragraph we recall the Mumford--Tate conjecture;
and in \cref{outline} we give an outline of the proof.
The ambitious reader may skip to \cref{Proof} and dive head first into the proof.
\begin{theorem}
	\label{MTCAxX}
	Let $K$ be a finitely generated subfield of~$\CC$.
	If $A$ is an abelian surface over~$K$ and $X$ is a K3~surface over~$K$,
	then the Mumford--Tate conjecture is true for $\HH^{2}(A \times X)(1)$.
\end{theorem}

\paragraph{The Mumford--Tate conjecture}
Let~$K$ be a finitely generated field of characteristic~$0$;
and let $\map{K}[hook]{\CC}$ be an embedding of~$K$ into the complex numbers.
Let~$\widebar{K}$ be the algebraic closure of~$K$ in~$\CC$.
Let~$X/K$ be a smooth projective variety.
One may attach several cohomology groups to~$X$.
For the purpose of this article we are interested in two cohomology theories:
Betti cohomology and $\ell$-adic \'{e}tale cohomology (for a prime number~$\ell$).
We will write $\HB^{i}(X)$ for the $\QQ$-Hodge structure $\HH_{\textnormal{sing}}^{i}(X(\CC), \QQ)$.
Similarly, we write $\Hl^{i}(X)$ for the $\Gal(\widebar{K}/K)$-representation $\HH_{\textnormal{\'{e}t}}^{i}(X_{\widebar{K}}, \QQl)$.

The Mumford--Tate conjecture is a precise way of saying that the cohomology
groups $\HB^{i}(X)$ and $\Hl^{i}(X)$ contain the same information about~$X$.
To make this precise, let $\GB(\HB^{i}(X))$ be the Mumford--Tate group of the Hodge structure~$\HB^{i}(X)$,
and let $\Glc(\Hl^{i}(X))$ be the connected component of the Zariski closure of $\Gal(\widebar{K}/K)$ in $\GL(\Hl^{i}(X))$.
The comparison theorem by Artin, comparing singular cohomology with \'{e}tale cohomology,
canonically identifies $\GL(\HB^{i}(X)) \otimes \QQl$ with $\GL(\Hl^{i}(X))$.
The Mumford--Tate conjecture (for the prime $\ell$, and the embedding $\map{K}[hook]{\CC}$) states that under this identification
\[
	\GB(\HB^{i}(X)) \otimes \QQl \cong \Glc(\Hl^{i}(X)).
\]

\paragraph{Outline of the proof}
\label{outline}
Let~$A/K$ be an abelian surface, and
let~$X/K$ be a K3~surface.
Observe that, by K\"{u}nneth's theorem, $\HB^{2}(A \times X) \cong \HB^{2}(A) \oplus \HB^{2}(X)$.
Similarly $\Hl^{2}(A \times X) \cong \Hl^{2}(A) \oplus \Hl^{2}(X)$.
Recall that the Mumford--Tate conjecture for~$A$ is known in degree~$1$, and hence in all degrees.
(This is classical, but see corollary~4.4 of~\cite{Lo14} for a reference.)
By~\cite{MTCK3I,MTCK3II,An96}, the Mumford--Tate conjecture for~$X$ (in degree~$2$) is true as well.
Still, it is not a formal consequence that the the Mumford--Tate conjecture for $A \times X$ is true in degree~$2$.

The proof of \cref{MTCAxX} falls apart into four cases, that use very different techniques.
All cases build on the Hodge theory of K3~surfaces and abelian varieties, of which we provide an overview in \cref{Hodge_theory}.

Let~$V$ be the transcendental part of~$\HB^{2}(X)$.
The first case (\cref{diffEnd}) inspects $\End(V)$,
and exploits Chebotaryov's density theorem, which we recall in \cref{End}.
The second case (\cref{diffLie}) looks at the Lie type of~$\GB(V)$,
and uses results about semisimple groups over number fields,
which we assemble in \cref{Lie}.

The third case (\cref{KSarg}) deals with Kummer varieties, and other K3~surfaces for which $\dim(V)$ is small.
We use the theory of Kuga--Satake varieties, and apply techniques of Lombardo, developed in~\cite{Lo14}.
The preliminaries of this part of the proof are gathered in \cref{AV}.

The final case (\cref{mtcnonsplitl}) is the only case where we use that $\HH^{2}(X)$ is a motive coming from a K3~surface.
We use information about the reduction of $X$ modulo a place of~$K$,
and combine this with a result about non-split groups
and results about compatible systems of $\ell$-adic representations.

\paragraph{Notation and terminology}
\label{intronot}
Let $K$ be a finitely generated field of characteristic~$0$;
and fix an embedding $\map{K}[hook]{\CC}$.
In this article we use the language of motives \`{a}~la Andr\'{e},~\cite{An95}.
To be precise, our category of base pieces is the category of smooth projective varieties over~$K$,
and our reference cohomology is Betti cohomology,~$\HB(\_)$; which, we stress, depends on the chosen embedding $\map{K}[hook]{\CC}$.
We write $\HH^{i}(X)$ for the motive of weight~$i$ associated with a smooth projective variety~$X/K$.

The Mumford--Tate conjecture naturally generalises to motives.
Let~$M$ be a motive.
We will write $\HB(M)$ for its Hodge realisation;
$\Hl(M)$~for its $\ell$-adic realisation;
$\GB(M)$~for its Mumford--Tate group (\textit{i.e.,} the Mumford--Tate group of~$\HB(M)$);
and $\Glc(M)$~for~$\Glc(\Hl(M))$.
We will use the notation $\MTC_{\ell}(M)$ for the conjectural statement
\[
	\GB(M) \otimes \QQl \cong \Glc(M),
\]
and $\MTC(M)$ for the assertion $\MTC_{\ell}(M)$ for all prime numbers~$\ell$.
In this paper, we never use specific properties of the chosen embedding $\map{K}[hook]{\CC}$,
and all statements are valid for every such embedding.
In particular, we will speak about subfields of~$\CC$, where the embedding is implicit.

In this paper, we will use compatible systems of $\ell$-adic representations.
We refer to the letters of Serre to Ribet (see~\cite{Serre}) or the work of Larsen and Pink \cite{LP2,LP3} for more information.

Throughout this paper, $A$~is an abelian variety, over some base field.
(Outside \cref{AV}, it is even an abelian surface.)
Assume $A$ is absolutely simple; and choose a polarisation of~$A$.
Let $(D, \dagger)$ be its endomorphism ring $\End^{0}(A)$
together with the Rosati involution associated with the polarisation.
The simple algebra~$D$ together with the positive involution~$\dagger$
has a certain type in the Albert classification that does not depend on the chosen polarisation.
We say that~$A$ is of type~\textsc{x} if $(D, \dagger)$ is of type~\textsc{x},
where~\textsc{x} runs through $\Set{\textsc{i}, \ldots, \textsc{iv}}$.
If $E$~denotes the center of~$D$, with degree $e = [E : \QQ]$,
we also say that $A$ is of type~$\textsc{x}(e)$.

Whenever we speak of (semi)simple groups or (semi)simple Lie algebras, we mean
\emph{non-commutative} (semi)simple groups, and \emph{non-abelian} (semi)simple
Lie-algebras.

Let~$T$ be a type of Dynkin diagram (\textit{e.g.,} $A_{n}$, $B_{n}$, $C_{n}$ or~$D_{n}$).
Let~$\mathfrak{g}$ be a semisimple Lie algebra over~$K$.
We say that $T$ does not occur in the Lie type of $\mathfrak{g}$,
if the Dynkin diagram of $\mathfrak{g}_{\widebar{K}}$ does not have a component of type~$T$.
For a semisimple group~$G$ over~$K$, we say that~$T$ does not occur in the Lie type of~$G$,
if~$T$ does not occur in the Lie type of $\Lie(G)$.

\paragraph{Acknowledgements}
I first and foremost want to thank Ben Moonen, my supervisor,
for his inspiration and help with critical parts of this paper.
Part of this work was done when the author was visiting Matteo Penegini at the
University of Milano; and I thank him for the hospitality and the inspiring collaboration.
I want to thank Bert van Geemen, Davide Lombardo,
and Milan Lopuha\"{a} for useful discussions about parts of the proof.
All my colleagues in Nijmegen who provided encouraging or insightful remarks
during the process of research and writing also deserve my thanks.
Further thanks goes to
grghxy, Guntram, and Mikhail Borovoi on \url{mathoverflow.net}%
\footnote{A preliminary version of \cref{trans} arose from a question on MathOverflow
	titled
	{``How simple does a $\mathbb{Q}$-simple group remain after base change
		to~$\mathbb{Q}_{\ell}$?''}
	(\url{http://mathoverflow.net/q/214603/78087}).
	The answers also inspired \cref{fixpoints}.}.

This research has been financially supported by the Netherlands Organisation for Scientific Research~(NWO) under 
project no.~613.001.207 \emph{(Arithmetic and motivic aspects of the Kuga--Satake construction)}.

\section{Some remarks on Chebotaryov's density theorem and transitive group actions}
\label{End}

\begin{theorem}[Chebotaryov's density theorem]
	\label{cheb}
	Let $K \subset E$ be an extension of number fields.
	Let $E \subset L$ be a Galois closure of~$E$, and
	let $G = \Gal(L/K)$ be the Galois group of~$L$ over~$K$.
	Let $\Sigma = \Hom_{K}(E, L)$ be the set of field embeddings over~$K$ of~$E$ in~$L$.
	\begin{itemize}
		\item Let $\mathfrak{p}$ be a prime of~$K$ that is unramified in~$L$, and
			let $C_{\mathfrak{p}} \subset G$ be the conjugacy class of the Frobenius elements associated with~$\mathfrak{p}$.
			The decomposition type of~$\mathfrak{p}$ in~$\mathcal{O}_{E}$
			is equal to the cycle type of~$C_{\mathfrak{p}}$ acting on~$\Sigma$.
		\item Let $C \subset G$ be a union of conjugacy classes of~$G$.
			The set
			\[
				\Set{ \mathfrak{p} \in \Spec(\mathcal{O}_{E}) \given
				\text{$\mathfrak{p}$ is unramified, and $C_{\mathfrak{p}} \subset C$} }
			\]
			has density~$\tfrac{\abs{C}}{\abs{G}}$ as subset of~$\Spec(\mathcal{O}_{E})$.
	\end{itemize}
	\begin{proof}
		See fact~2.1 and theorem~3.1 of~\cite{Cheb}.
		See Theorem~13.4 of~\cite{Neu} for the case where~$E/K$ is Galois.
	\end{proof}
\end{theorem}

\begin{lemma}
	\label{fixpoints}
	Let $G$ be a finite group acting transitively on a finite set~$\Sigma$.
	Let~$n \in \ZZ_{\ge0}$ be a non-negative integer, and
	let $C \subset G$ be the set of elements $g \in G$ that have at least~$n$ fixed points:
	\[
		C = \Set[\big]{ g \in G \given \abs{\Sigma^{g}} \ge n }
	\]
	If $n \cdot \abs{C} \ge \abs{G}$, then $\abs{\Sigma} = n$.
	If furthermore the action of~$G$ on~$\Sigma$ is faithful, then $\abs{G} = n$, and~$\Sigma$ is principal homogeneous under~$G$.
	\begin{proof}
		Burnside's lemma gives
		\[
			1 = \abs{G\backslash \Sigma} =
			\frac{1}{\abs{G}} \sum_{g \in G} \abs{\Sigma^{g}} \ge \frac{n \cdot \abs{C}}{\abs{G}} \ge 1.
		\]
		Hence $n \cdot \abs{C} = \abs{G}$ and all elements in~$C$ have exactly~$n$ fixed points.
		In particular the identity element has~$n$ fixed points, which implies $\abs{\Sigma} = n$.
		If $G$ acts faithfully on~$\Sigma$, then $\abs{\Sigma} = n$ implies $C = \Set{e}$,
		and thus $\abs{G} = n = \abs{\Sigma}$. So~$\Sigma$ is principal homogeneous under~$G$.
	\end{proof}
\end{lemma}

\begin{lemma}
	\label{trans}
	Let $F_{1}$ be a Galois extension of $\QQ$.
	Let $F_{2}$ be a number field.
	If for all prime numbers~$\ell$,
	the product of local fields $F_{1} \otimes \QQl$ is a factor of $F_{2} \otimes \QQl$,
	then $F_{1} \cong F_{2}$.
	\begin{proof}
		Let $L$ be a Galois closure of~$F_{2}$, and
		let $G$ be the Galois group~$\Gal(L/\QQ)$,
		which acts naturally on the set of field embeddings $\Sigma = \Hom(F_{2}, L)$.
		Let~$n$ be the degree of~$F_{1}$, and
		let~$C$ be the set $\Set[\big]{ g \in G \given \abs{\Sigma^{g}} \ge n }$
		of elements in~$G$ that have at least~$n$ fixed points in~$\Sigma$.

		By Chebotaryov's density theorem (\labelcref{cheb}),
		the set of primes that split completely in~$F_{1}$ has density~$1/n$.
		Another application of \cref{cheb} shows that the set of primes~$\ell$
		for which $F_{2} \otimes \QQl$ has a semisimple factor isomorphic to~$(\QQl)^{n}$ must have density~$\ge 1/n$.
		Our assumption therefore implies that $n \cdot \abs{C} \ge \abs{G}$.
		By \cref{fixpoints}, this implies $\abs{\Sigma} = n$, and since~$G$ acts faithfully on~$\Sigma$,
		we find that $F_{2}/\QQ$ is Galois of degree~$n$.
		Because Galois extensions of number fields can be recovered from their set of splitting primes
		(Satz~VII.13.9 of~\cite{Neu}), we conclude that $F_{2} \cong F_{1}$.
	\end{proof}
\end{lemma}

\begin{lemma}
	\label{hack}
	Let $F_{1}$ be a quadratic extension of~$\QQ$.
	Let $F_{2}$ be a number field of degree~$\le 5$ over~$\QQ$.
	If for all prime numbers~$\ell$, the products of local fields
	$F_{1} \otimes \QQl$ and $F_{2} \otimes \QQl$ have an isomorphic factor,
	then $F_{1} \cong F_{2}$.
	\begin{proof}
		Let $L$ be a Galois closure of~$F_{2}$, and
		let $G$ be the Galois group~$\Gal(L/\QQ)$,
		which acts naturally on the set of field embeddings $\Sigma = \Hom(F_{2}, L)$.
		Observe that $G$ acts transitively on~$\Sigma$,
		and we identify~$G$ with its image in~$\mathfrak{S}(\Sigma)$.
		Write $n$ for the degree of~$F_{2}$ over~$\QQ$, which also equals~$\abs{\Sigma}$.
		The order of~$G$ is divisible by~$n$.
		Hence, if~$n$ is prime, then $G$ must contain an $n$-cycle.
		
		Suppose that $G$ contains an $n$-cycle.
		By Chebotaryov's density theorem (\labelcref{cheb})
		there must be a prime number~$\ell$ that is inert in~$F_{2}$.
		By our assumption $F_{2} \otimes \QQl$ also contains a factor of at most degree~$2$ over~$\QQl$.
		This shows that $n = 2$.

		If $n = 4$, then $G$ does not contain an $n$-cycle if and only if it is isomorphic to~$V_{4}$ or~$A_{4}$.
		If $G \cong V_{4}$, only the identity element has fixed points,
		and by Chebotaryov's density theorem this means that
		the set of primes~$\ell$ for which $F_{2} \otimes \QQl$ has a factor~$\QQl$ has density~$1/4$,
		whereas the set of primes splitting in~$F_{1}$ has density~$1/2$.
		On the other hand,
		if $G \cong A_{4}$, only~$3$ of the~$12$ elements have a $2$-cycle in the cycle decomposition,
		and by Chebotaryov's density theorem this means that
		the set of primes~$\ell$ for which $F_{2} \otimes \QQl$
		has a factor isomorphic to a quadratic extension of~$\QQl$ has density~$1/4$,
		whereas the set of primes inert in~$F_{1}$ has density~$1/2$.
		This gives a contradiction.
		We conclude that~$n$ must be~$2$; and therefore $F_{1} \cong F_{2}$, by \cref{trans}.
	\end{proof}
\end{lemma}

\section{Several results on semisimple groups over number fields}
\label{Lie}

Throughout this section $K$ is a field of characteristic~$0$.
\begin{lemma}
	\label{isomLie}
	Let $G$ be a connected algebraic group over~$K$, and
	let $H \subset G$ be a subgroup.
	If $\Lie(H) = \Lie(G)$, then $H = G$.
	\begin{proof}
		This is immediate, since $H$ is a subgroup of~$G$ of the same dimension as~$G$.
	\end{proof}
\end{lemma}

\begin{lemma}[Goursat's lemma for Lie algebras]
	\label{Goursat}
	Let $\mathfrak{g}_{1}$ and $\mathfrak{g}_{2}$ be Lie algebras over~$K$, and
	let $\mathfrak{h} \subset \mathfrak{g}_{1} \oplus \mathfrak{g}_{2}$ be a sub-Lie algebra
	such that the projections $\map[\pi_{i}]{\mathfrak{h}}{\mathfrak{g}_{i}}$ are surjective.
	Let $\mathfrak{n}_{1}$ be the kernel of~$\pi_{2}$, and~$\mathfrak{n}_{2}$ the kernel of~$\pi_{1}$.
	The projection~$\pi_{i}$ identifies~$\mathfrak{n}_{i}$ with an ideal of~$\mathfrak{g}_{i}$,
	and the image of the canonical map
	\[
		\map{\mathfrak{h}}
		{(\mathfrak{g}_{1}/\pi_{1}(\mathfrak{n}_{1})) \oplus (\mathfrak{g}_{2}/\pi_{2}(\mathfrak{n}_{2}))}
	\]
	is the graph of an isomorphism
	$\map{\mathfrak{g}_{1}/\pi_{1}(\mathfrak{n}_{1})}{\mathfrak{g}_{2}/\pi_{2}(\mathfrak{n}_{2})}$.
	\begin{proof}
		Observe that~$\pi_{i}$ is injective on~$\mathfrak{n}_{i}$.
		If $x \in \pi_{i}(\mathfrak{n}_{i})$ and $y \in \mathfrak{g}_{i}$,
		then $[x,y] \in \pi_{i}(\mathfrak{n}_{i}\mkern-1mu)$, because $\pi_{i}$ is surjective,
		and $\mathfrak{n}_{i}$ is an ideal of~$\mathfrak{h}$.
		Let $\widebar{\mathfrak{h}}$ be the image of the canonical map
		\[
			\map{\mathfrak{h}}
			{(\mathfrak{g}_{1}/\pi_{1}(\mathfrak{n}_{1})) \oplus (\mathfrak{g}_{2}/\pi_{2}(\mathfrak{n}_{2}))}
		\]
		By construction, the projections $\map{\widebar{\mathfrak{h}}}{\mathfrak{g}_{i}/\pi_{i}(\mathfrak{n}_{i})}$ are injective;
		and they are surjective by assumption.
		This proves the lemma.
	\end{proof}
\end{lemma}

\begin{remark}
	\label{summand}
	Let $\mathfrak{h} \subset \mathfrak{g}_{1} \oplus \mathfrak{g}_{2}$
	be Lie algebras over~$K$ satisfying the conditions of \cref{Goursat}.
	Assume that $\mathfrak{g}_{1}$ and~$\mathfrak{g}_{2}$ are finite-dimensional and semisimple.
	It follows from the proof of \cref{Goursat} that
	there exist semisimple Lie algebras $\mathfrak{s}_{1}$, $\mathfrak{t}$, and~$\mathfrak{s}_{2}$
	such that $\mathfrak{g}_{1} \cong \mathfrak{s}_{1} \oplus \mathfrak{t}$,
	$\mathfrak{g}_{2} \cong \mathfrak{t} \oplus \mathfrak{s}_{2}$,
	and $\mathfrak{h} \cong \mathfrak{s}_{1} \oplus \mathfrak{t} \oplus \mathfrak{s}_{2}$.
\end{remark}

\begin{corollary}
	\label{Lie_obs}
	Let $K \subset L$ be a field extension.
	Let $G_{1}$ and~$G_{2}$ be connected semisimple groups over~$K$.
	Let $\map[\iota]{G}[hook]{G_{1} \times G_{2}}$ be a subgroup,
	with surjective projections onto both factors.
	If $\Lie(G_{1})_{L}$ and $\Lie(G_{2})_{L}$ have no isomorphic factor over~$L$,
	then $\iota$ is an isomorphism.
\end{corollary}

\begin{lemma}
	\label{normLie}
	Let $K \subset F$ be a finite field extension.
	Let $G$ be an algebraic group over~$F$.
	The Lie algebra $\Lie(\Res_{F/K}G)$ is isomorphic to the Lie algebra $\Lie(G)$, viewed as Lie algebra over~$K$.
	\begin{proof}
		This follows from the following diagram, the rows of which are exact.
		\[
			\begin{tikzcd}
				0 \rar & \Lie(\Res_{F/K}G) \rar & (\Res_{F/K}G)(K[\varepsilon]) \rar \dar{\text{\rotatebox{-90}{$\simeq$}}}
				& (\Res_{F/K}G)(K) \rar \dar{\text{\rotatebox{-90}{$\simeq$}}} & 0 \\
				0 \rar & \Lie(G) \rar & G(K[\varepsilon]) \rar & G(K) \rar & 0
			\end{tikzcd}
		\]
	\end{proof}
\end{lemma}

\begin{lemma}
	\label{Galois_obs}
	Let $F_{1}/K$ and $F_{2}/K$ be finite field extensions.
	Let $\mathfrak{g}_{i}/F_{i}$ ($i = 1,2$) be a finite product of absolutely simple Lie algebras
	(\textit{cf.}\ our conventions in \cref{intronot}).
	Write $(\mathfrak{g}_{i})_{K}$ for the Lie algebra~$\mathfrak{g}_{i}$ viewed as Lie algebra over~$K$.
	If $(\mathfrak{g}_{1})_{K}$ and~$(\mathfrak{g}_{2})_{K}$ have an isomorphic factor,
	then $F_{1} \cong_{K} F_{2}$.
	\begin{proof}
		The $K$-simple factors of $(\mathfrak{g}_{i})_{K}$ are all of the form $(\mathfrak{t}_{i})_{K}$,
		where $\mathfrak{t}_{i}$~is an $F_{i}$-simple factor of~$\mathfrak{g}_{i}$.
		So if $(\mathfrak{g}_{1})_{K}$ and~$(\mathfrak{g}_{2})_{K}$ have an isomorphic factor,
		there exist $F_{i}$-simple factors~$\mathfrak{t}_{i}$ of~$\mathfrak{g}_{i}$
		for which there exists an isomorphism $\map[f]{(\mathfrak{t}_{1})_{K}}{(\mathfrak{t}_{2})_{K}}$.
		Let $\widebar{K}$ be an algebraic closure of~$K$.
		Observe that
		\[
			(\mathfrak{t}_{i})_{K} \otimes_{K} \widebar{K} \cong
			\bigoplus_{\sigma \in \Hom_{K}(F_{i}, \widebar{K})} \mathfrak{t}_{i} \otimes_{F_{i},\sigma} \widebar{K},
		\]
		and note that $\Gal(\widebar{K}/K)$ acts transitively on $\Hom_{K}(F_{i}, \widebar{K})$.
		By assumption, the $\mathfrak{t}_{i}$ are $F_{i}$-simple,
		and therefore the $\mathfrak{t}_{i} \otimes_{F_{i},\sigma} \widebar{K}$
		are precisely the simple ideals of~$(\mathfrak{t}_{i})_{K} \otimes_{K} \widebar{K}$.
		Thus the isomorphism $f$ gives a $\Gal(\widebar{K}/K)$-equivariant bijection between the simple ideals
		of~$(\mathfrak{t}_{1})_{K} \otimes_{K} \widebar{K}$ and~$(\mathfrak{t}_{2})_{K} \otimes_{K} \widebar{K}$;
		and therefore $\Hom_{K}(F_{1}, \widebar{K})$ and~$\Hom_{K}(F_{2}, \widebar{K})$
		are isomorphic as $\Gal(\widebar{K}/K)$-sets.
		This proves the result.
	\end{proof}
\end{lemma}

\begin{lemma}
	\label{Galois_obs_ell}
	Let $F_{1}$ and $F_{2}$ be number fields.
	Let $G_{i}/F_{i}$ ($i = 1,2$) be an almost direct product of connected absolutely simple $F_{i}$-groups.
	Let~$\ell$ be a prime number, and
	let $\map[\iota_{\ell}]{G}[hook]{(\Res_{F_{1}/\QQ}G_{1})_{\QQl} \times (\Res_{F_{2}/\QQ}G_{2})_{\QQl}}$
	be a subgroup over~$\QQl$, with surjective projections onto both factors.
	If $\iota_{\ell}$ is not an isomorphism, then $F_{1} \otimes \QQl$ and $F_{2} \otimes \QQl$ have an isomorphic simple factor.
	\begin{proof}
		Observe that $(\Res_{F_{i}/\QQ}G_{i}) \otimes \QQl \cong
		\prod_{\lambda\mid\ell} \Res_{F_{i,\lambda}/\QQl} (G_{i} \otimes_{F_{i}} F_{\lambda})$.
		If $\iota_{\ell}$ is not an isomorphism, then by \cref{Lie_obs}, there exist places~$\lambda_{i}$ of~$F_{i}$ over~$\ell$
		such that $\Lie(\Res_{F_{1,\lambda_{1}}/\QQl}(G_{1} \otimes_{F_{1}} F_{1,\lambda_{1}}))$ and
		$\Lie(\Res_{F_{2,\lambda_{2}}/\QQl}(G_{2} \otimes_{F_{2}} F_{2,\lambda_{2}}))$
		have an isomorphic factor.
		By \cref{normLie,Galois_obs}, this implies that $F_{1,\lambda_{1}} \cong_{\QQl} F_{2,\lambda_{2}}$,
		which proves the lemma.
	\end{proof}
\end{lemma}

\section{Several results on abelian motives}
\label{AbMot}

\begin{lemma}
	\label{MTCcentres}
	The Mumford--Tate conjecture on centres is true for abelian motives.
	In other words, let $M$ be an abelian motive.
	Let $Z_{\textnormal{B}}(M)$ be the centre of the Mumford--Tate group~$\GB(M)$, and
	let $Z_{\ell}(M)$ be the centre of~$\Glc(M)$.
	Then $Z_{\ell}(M) \cong Z_{\textnormal{B}}(M) \otimes \mathbb{Q}_{\ell}$.
	\begin{proof}
		The result is true for abelian varieties
		(see theorem~1.3.1 of~\cite{Va08} or corollary~2.11 of~\cite{UY13}).

		By definition of abelian motive,
		there is an abelian variety~$A$ such that
		$M$ is contained in the Tannakian subcategory of motives generated by~$\HH(A)$.
		This yields a surjection $\map{\GB(A)}[twohead]{\GB(M)}$,
		and therefore $Z_{\textnormal{B}}(M)$ is the image of $Z_{\textnormal{B}}(A)$ under this map.
		The same is true on the $\ell$-adic side.
		Thus we obtain a commutative diagram with solid arrows
		\[
			\begin{tikzcd}
				Z_{\ell}(A) \dar{\text{\rotatebox{-90}{$\simeq$}}} \rar[two heads] &
				Z_{\ell}(M) \dar[dotted]{\text{\rotatebox{-90}{$\simeq$}}} \rar[hook] &
				\Glc(M) \dar[hook] \\
				Z_{\textnormal{B}}(A) \otimes \QQl \rar[two heads] &
				Z_{\textnormal{B}}(M) \otimes \QQl \rar[hook] &
				\GB(M) \otimes \QQl
			\end{tikzcd}
		\]
		which shows that the dotted arrow exists and is an isomorphism.
	\end{proof}
\end{lemma}

\begin{lemma}
	\label{MTC_Tate_indep}
	Let $\mathbb{1}$ denote the trivial motive.
	If $M$ is a motive, then the Mumford--Tate conjecture for~$M$
	is equivalent to the Mumford--Tate conjecture for~$M \oplus \mathbb{1}$.
	\begin{proof}
		Indeed, $M$ and $M \oplus \mathbb{1}$ generate the same Tannakian subcategory of motives.
	\end{proof}
\end{lemma}

\begin{lemma}
	\label{fgext}
	Let $K \subset L$ be an extension of finitely generated subfields of~$\CC$.
	If $M$ is a motive over~$K$, then $\MTC(M) \iff \MTC(M_{L})$.
	\begin{proof}
		See proposition~1.3 of~\cite{Mo15}.
	\end{proof}
\end{lemma}

\begin{lemma}
	\label{MTC_lindep}
	Let $M$ be an abelian motive.
	Assume that the $\ell$-adic realisations of~$M$ form a compatible system of $\ell$-adic representations.
	If the Mumford--Tate conjecture for~$M$ is true for one prime~$\ell'$, then it is true for all primes~$\ell$.
	\begin{proof}
		Since $M$ is an abelian motive, we have $\Glc(M) \subset \GB(M) \otimes \mathbb{Q}_{\ell}$.
		By our assumption on the $\ell$-adic realisations of~$M$,
		the proofs of theorem~4.3 and lemma~4.4 of~\cite{LP3} apply verbatim to our situation.
	\end{proof}
\end{lemma}

\paragraph{}
\label[conditioni]{star}
Let $K$ be a finitely generated subfield of~$\CC$.
A pair $(A,X)$, consisting of an abelian surface~$A$ and a K3~surface~$X$ over~$K$,
is said to satisfy \cref{star} for~$\ell$ if
\[
	\map{\Glc\big(\HH^{2}(A \times X)(1)\big)^{\der}}[hook]
	{\Glc\big(\HH^{2}(A)(1)\big)^{\der} \times \Glc\big(\HH^{2}(X)(1)\big)^{\der}}
\]
is an isomorphism.
\begin{lemma}
	\label{cata}
	Let $K$ be a finitely generated subfield of~$\CC$, and
	let $C/K$ be a smooth (not necessarily proper) curve over~$K$, with generic point~$\eta$.
	Let $A/C$ be an abelian scheme, and let $X/C$ be a K3~surface.
	There exists a closed point $c \in C$ and a prime number~$\ell$ such that
	$\Glc(\HH^{2}(A_{\eta})(1)) \cong \Glc(\HH^{2}(A_{c})(1))$ and
	$\Glc(\HH^{2}(X_{\eta})(1)) \cong \Glc(\HH^{2}(X_{c})(1))$.
	Furthermore, if $(A_{c},X_{c})$ satisfies \cref{star} for~$\ell$,
	then so does~$(A_{\eta},X_{\eta})$.
	\begin{proof}
		The existence of the point~$c$ follows immediately from theorem~1.1 of~\cite{CaTa}.
		The diagram
		\[
			\begin{tikzcd}
				\Glc(\HH^{2}(A_{c} \times X_{c})(1))^{\der} \rar[hook] \dar[hook] &
				\Glc(\HH^{2}(A_{\eta} \times X_{\eta})(1))^{\der} \dar[hook] \\
				\Glc(\HH^{2}(A_{c})(1))^{\der} \times \Glc(\HH^{2}(X_{c})(1))^{\der}\mkern-1.5mu \rar{\simeq} &
				\Glc(\HH^{2}(A_{\eta})(1))^{\der} \times \Glc(\HH^{2}(X_{\eta})(1))^{\der}
			\end{tikzcd}
		\]
		shows that $(A_{\eta},X_{\eta})$ satisfies \cref{star} for~$\ell$ if $(A_{c},X_{c})$ satisfies it.
	\end{proof}
\end{lemma}

\paragraph{}
Let $\ell$ be a prime number.
Let $G_{1}$ and $G_{2}$ be connected reductive groups over~$\QQl$.
By a \emph{$(G_{1},G_{2})$-tuple} over~$K$ we shall mean a pair~$(A,X)$,
where $A$ is an abelian surface over~$K$, and~$X$ is a K3~surface over~$K$
such that $\Glc(\HH^{2}(A)(1)) \cong G_{1}$ and $\Glc(\HH^{2}(X)(1)) \cong G_{2}$.
We will show in \cref{Proof} that there exist groups~$G_{1}$ and~$G_{2}$
that satisfy the hypothesis of the following lemma,
namely that \cref{star} for~$\ell$ is satisfied for all $(G_{1},G_{2})$-tuples over number fields.

\begin{lemma}
	\label{specialise}
	Let $\ell$ be a prime number.
	Let $G_{1}$ and~$G_{2}$ be connected reductive groups over~$\QQl$.
	If for all number fields~$K$, all $(G_{1},G_{2})$-tuples~$(A,X)$ over~$K$
	satisfy \cref{star} for~$\ell$,
	then for all finitely generated subfields~$L$ of~$\CC$, all $(G_{1},G_{2})$-tuples~$(A,X)$ over~$L$
	satisfy \cref{star} for~$\ell$.
	\begin{proof}
		The proof goes by induction on the transcendence degree~$n$ of~$L$.
		If $n = 0$, the result is true by assumption.
		Suppose that $n > 0$,
		and assume as induction hypothesis that \cref{star} for~$\ell$ is satisfied for
		all $(G_{1},G_{2})$-tuples over all finitely generated subfields of~$\CC$
		with transcendence degree~$< n$.

		There exists a field $K \subset L$, and a smooth curve~$C/K$ such that~$L$ is the function field of~$C$.
		Observe that $\trdeg(K) = n - 1$.
		By the induction hypothesis and \cref{cata}, the claim of the lemma is true for~$L$.
		The result follows by induction.
	\end{proof}
\end{lemma}

\section{Some remarks on the Mumford--Tate conjecture for abelian varieties}
\label{AV}

\paragraph{} For the convenience of the reader, we copy some results from \cite{Lo14}.
Before we do that, let us recall the notion of the \emph{Hodge group}, $\HdgB(A)$, of an abelian variety.
Let $A$ be an abelian variety over a finitely generated field $K \subset \CC$.
By definition, the Mumford--Tate group of an abelian variety is $\GB(A) = \GB(\HB^{1}(A)) \subset \GL(\HB^{1}(A))$, and
we put
\[
	\HdgB(A) = (\GB(A) \cap \SL(\HB^{1}(A)))^{\circ} \quad \text{and} \quad
	\Hdgl(A) = (\Gl(A) \cap \SL(\Hl^{1}(A)))^{\circ}.
\]
We leave it as an easy exercise to the reader to verify that
\[
	\MTC_{\ell}(A) \quad\iff\quad \HdgB(A) \otimes \QQl \cong \Hdgl(A).
\]

\begin{definition}[1.1 in \cite{Lo14}]
	\label{reldim}
	Let $A$ be an absolutely simple abelian variety of dimension~$g$ over~$K$.
	The endomorphism ring $D = \End^{0}(A)$ is a division algebra.
	Write $E$ for the centre of~$D$.
	The ring $E$ is a field, either \tr{} (totally real) or~\cm.
	Write $e$ for $[E : \QQ]$.
	The degree of~$D$ over~$E$ is a perfect square~$d^{2}$.

	The \emph{relative dimension} of~$A$ is
	\[
		\reldim(A) =
		\begin{cases}
			\frac{g}{de}, &\text{if $A$ is of type \textsc{i}, \textsc{ii}, or \textsc{iii}}, \\
			\frac{2g}{de}, &\text{if $A$ is of type \textsc{iv}}. \\
		\end{cases}
	\]
\end{definition}
Note that $d = 1$ if $A$ is of type \textsc{i}, and~$d = 2$ if $A$ is of type \textsc{ii} or~\textsc{iii}.

In definition~2.22 of \cite{Lo14}, Lombardo defines when an abelian variety is of \emph{general Lefschetz type}.
This definition is a bit unwieldy, and its details do not matter too much for our purposes.
What matters are the following results, that prove that certain abelian varieties are of general Lefschetz type,
and that show why this notion is relevant for us.

\begin{lemma} \label{gLtAVs}
	Let $A$ be an absolutely simple abelian variety over a finitely generated subfield of~$\CC$.
	Assume that $A$ is of type \textsc{i} or~\textsc{ii}.
	If $\reldim(A)$ is odd, or equal to~$2$, then $A$ is of general Lefschetz type.
	\begin{proof}
		If $\reldim(A)$ is odd, then this follows from theorems~6.9 and~7.12 of~\cite{BGK}.
		Lombardo notes (remark~2.25 in~\cite{Lo14})
		that the proof of \cite{BGK} also works if $\reldim(A) = 2$,
		and also refers to theorem~8.5 of~\cite{Ch92} for a proof of that fact.
	\end{proof}
\end{lemma}

\begin{lemma}
	Let $K$ be a finitely generated subfield of~$\CC$.
	Let $A_{1}$ and~$A_{2}$ be two abelian varieties over~$K$ that are
	isogenous to products of abelian varieties of general Lefschetz type.
	If $D_{4}$ does not occur in the Lie type of $\Hdgl(A_{1})$ and~$\Hdgl(A_{2})$,
	then either
	\[
		\Hom_{K}(A_{1}, A_{2}) \ne 0,
		\quad \text{or} \quad
		\Hdgl(A_{1} \times A_{2}) \cong \Hdgl(A_{1}) \times \Hdgl(A_{2}).
	\]
	\begin{proof}
		This is remark~4.3 of~\cite{Lo14},
		where Lombardo observes that, under the assumption of the
		lemma, theorem~4.1 of~\cite{Lo14} can be applied to products of
		abelian varieties of general Lefschetz type.
	\end{proof}
\end{lemma}

\begin{lemma} \label{MTCgLt}
	Let $A$ be an abelian variety over a finitely generated field $K \subset \CC$.
	Let $L \subset \CC$ be a finite extension of~$K$
	for which $A_{L}$ is isogenous over~$L$ to a product of absolutely simple abelian varieties $\prod A_{i}^{k_{i}}$.
	Assume that for all~$i$ the following conditions are valid:
	\begin{condition}[label=(\alph*)]
	\item either $A_{i}$ is of general Lefschetz type or $A_{i}$ is of \cm{}~type;
	\item the Lie type of $\Hdgl(A_{i})$ does not contain~$D_{4}$;
	\item \label{MTCfactors}
		the Mumford--Tate conjecture is true for~$A_{i}$.
	\end{condition}
	Under these conditions the Mumford--Tate conjecture is true for~$A$.
	\begin{proof}
		By \cref{fgext} we know that $\MTC(A) \iff \MTC(A_{L})$.
		Furthermore, note that $\MTC(A_{L})$ is equivalent to~$\MTC(\prod A_{i})$.
		Observe that
		\[
			\Hdgl(A) \subset \HdgB(A) \otimes \QQ_{\ell} \subset \prod \HdgB(A_{i}) \otimes \QQ_{\ell} = \prod \Hdgl(A_{i}),
		\]
		where the first inclusion is
		Deligne's ``Hodge = absolute Hodge'' theorem;
		the second inclusion is a generality;
		and the last equality is \cref{MTCfactors}.

		If we ignore the factors that are~\cm, then an inductive
		application of the previous lemma yields $\Hdgl(A) = \prod \Hdgl(A_{i})$.
		If we do not ignore the factors that are~\cm,
		then we actually get $\Hdgl(A)^{\der} = \prod \Hdgl(A_{i})^{\der}$.
		Together with \cref{MTCcentres},
		this proves $\Hdgl(A) = \HdgB(A) \otimes \QQ_{\ell}$.
	\end{proof}
\end{lemma}
As an illustrative application of this result, Lombardo observes in
corollary~4.5 of~\cite{Lo14} that the Mumford--Tate conjecture is true for
arbitrary products of elliptic curves and abelian surfaces.

\section{Hodge theory of K3~surfaces and abelian surfaces}
\label{Hodge_theory}

\paragraph{}
In this section we recall some results of Zarhin
that describe all possible Mumford--Tate groups of Hodge structures of K3~type,
\textit{i.e.,} Hodge structures of weight~$0$ with Hodge numbers of the form $(1,n,1)$.

The canonical example of a Hodge structure of K3~type is the cohomology in degree~$2$ of a complex K3~surface~$X$.
Namely the Hodge structure~$\HB^{2}(X)(1)$ has Hodge numbers $(1,20,1)$.
Another example is provided by abelian surfaces, which is the content of \cref{Hodge_theory_A} below.

\pagebreak[2]
\begin{lemma}
	\label{Hodge_theory_K3}
	Let $V$ be an irreducible Hodge structure of K3~type, and
	let $\psi$ be a polarisation on~$V$.
	\begin{enumerate}
		\item The endomorphism algebra~$E$ of~$V$ is a field.
		\item The field~$E$ is~\tr{} (totally real) or~\cm.
		\item If $E$ is~\tr{}, then $\dim_{E}(V) \ge 3$.
		\item Let $\tilde{\psi}$ be the unique $E$-bilinear (resp.\ hermitian) form 
			such that $\psi = \trace_{E/\QQ} \circ \tilde{\psi}$ if~$E$ is~\tr{} (resp.~\cm).
			Let $E_{0}$ be the maximal totally real subfield of~$E$.
			The Mumford--Tate group of~$V$ is
			\[
				\GB(V) \cong
				\begin{cases}
					\Res_{E/\QQ}\SO(\tilde{\psi}), &\text{if $E$ is \tr;} \\
					\Res_{E_{0}/\QQ}\U_{E/E_{0}}(\tilde{\psi}), &\text{if $E$ is \cm.}
				\end{cases}
			\]
	\end{enumerate}
	\begin{proof}
		The first (resp.\ second) claim is theorem~1.6.a (resp.\ theorem~1.5) of~\cite{Za83};
		the third claim is observed by Van~Geemen, in lemma~3.2 of \cite{vG08}; and
		the final claim is a combination of theorems~2.2 and~2.3 of~\cite{Za83}.
		(We note that \cite{Za83} deals with Hodge groups,
		but because our Hodge structure has weight~$0$,
		the Mumford--Tate group and the Hodge group coincide.)
	\end{proof}
\end{lemma}

\begin{remark}
	\label{SO4}
	Let $V$, $E$ and~$\tilde{\psi}$ be as in \cref{Hodge_theory_K3}.
	If $E$ is~\cm, then $\U(\tilde{\psi})^{\der} = \SU(\tilde{\psi})$ is absolutely simple over~$E$.
	If $E$ is~\tr{} and $\dim_{E}(V) \ne 4$, then $\SO(\tilde{\psi})$ is absolutely simple over~$E$.
	Assume $E$ is~\tr{} and $\dim_{E}(V) = 4$.
	In this case $\SO(\tilde{\psi})$ is not absolutely simple over~$E$;
	it has Lie type $D_{2} = A_{1} \oplus A_{1}$.
	In this remark we will take a close look at this special case,
	because a good understanding of it will play a crucial r\^{o}le in the proof of \cref{mtcnonsplitl}.

	Geometrically we find
	$\SO(\tilde{\psi})_{\widebar{E}} \cong (\SL_{2,\widebar{E}} \times \SL_{2,\widebar{E}})/ \langle (-1,-1) \rangle$.
	We distinguish the following two cases:
	\begin{case}
		\item \label{SO4nonsimple} $\SO(\tilde{\psi})$ is not simple over~$E$.
			The fact that is most relevant to us is that there exists a quaternion algebra~$D/E$
			such that $\SO(\tilde{\psi}) \cong (N \times N^{\textnormal{op}})/\langle (-1,-1) \rangle$
			where $N$ is the group over~$E$ of elements in $D^{\star}$ that have norm~$1$,
			and $N^{\textnormal{op}} \subset (D^{\textnormal{op}})^{\star}$
			is the group of units with norm~$1$ in $D^{\textnormal{op}}$.
			One can read more about the details of this claim in section~8.1 of~\cite{Mo15}.
			This situation is also described in section~26.B of~\cite{BoI},
			where the quaternion algebra is replaced by $D \times D$ viewed as quaternion algebra over $E \times E$.
			This might be slightly more natural,
			but it requires bookkeeping of \'{e}tale algebras which
			makes the proof in \cref{Proof} more difficult than necessary.
		\item \label{SO4simple} $\SO(\tilde{\psi})$ is simple over~$E$.
			This means that the action of~$\Gal(\widebar{E}/E)$ on~$\SO(\tilde{\psi})_{\widebar{E}}$
			interchanges the two factors~$\SL_{2,\widebar{E}}$.
			The stabilisers of these factors are subgroups of index~$2$ that coincide.
			This subgroup fixes a quadratic extension~$F/E$.
			From our description of the geometric situation,
			together with the description of the stabilisers,
			we see that $\textnormal{Spin}(\tilde{\psi}) = \Res_{F/E}\scrG$
			is a $(2 : 1)$-cover of~$\SO(\tilde{\psi})$,
			where $\scrG$~is an absolutely simple, simply connected group of Lie type~$A_{1}$ over~$F$.
	\end{case}
	What we have gained is that in all cases we have a description (up to isogeny) of~$\GB(V)^{\der}$
	as Weil restriction of a group that is an \emph{almost direct product} of groups that are \emph{absolutely simple}.
	This allows us to apply \cref{Galois_obs}, which will play an important r\^{o}le in \cref{Proof}.
\end{remark}

\paragraph{Notation and terminology}
\label{FG}
Let $V$, $E$ and~$\tilde{\psi}$ be as in \cref{Hodge_theory_K3}.
To harmonise the proof in \cref{Proof}, we unify notation as follows:
\[
	F = 
	\begin{cases}
		E_{0} & \text{if $E$ is \cm,} \\
		E & \text{if $E$ is \tr{} and $\dim_{E}(V) \ne 4$,} \\
		E & \text{if $E$ is \tr, $\dim_{E}(V) = 4$, and we are in \cref{SO4nonsimple},} \\
		F & \text{if $E$ is \tr, $\dim_{E}(V) = 4$, and we are in \cref{SO4simple}.}
	\end{cases}
\]
Similarly
\[
	\scrG = 
	\begin{cases}
		\U(\tilde{\psi}) & \text{if $E$ is \cm,} \\
		\SO(\tilde{\psi}) & \text{if $E$ is \tr{} and $\dim_{E}(V) \ne 4$,} \\
		\SO(\tilde{\psi}) & \text{if $E$ is \tr, $\dim_{E}(V) = 4$, and we are in \cref{SO4nonsimple},} \\
		\scrG & \text{if $E$ is \tr, $\dim_{E}(V) = 4$, and we are in \cref{SO4simple}.}
	\end{cases}
\]
We stress that $\scrG^{\der}$ is an almost direct product of absolutely simple groups over~$F$.
In \cref{Proof}, most of the time it is enough to know that $\GB(V)$~is isogenous to~$\Res_{F/\QQ}\scrG$.
When we need more detailed information, it is precisely the case that $E$~is \tr{} and $\dim_{E}(V) = 4$.
For this case we gave a description of~$\scrG$ in the previous remark.

\paragraph{}
Let $V$, $E$ and~$\tilde{\psi}$ be as in \cref{Hodge_theory_K3}.
Write $n$ for~$\dim_{E}(V)$.
If $E$ is~\tr, then we say that the group~$\SO(\tilde{\psi})$ over~$E$ is a group of type~$\SO_{n,E}$.
We also say that $\GB(V)$ is of type~$\Res_{E/\QQ}\SO_{n,E}$.
Similarly, if $E$ is~\cm{}, with maximal totally real subfield~$E_{0}$,
then we say that the group~$\U_{E/E_{0}}(\tilde{\psi})$ over~$E_{0}$ is a group of type~$\U_{n,E_{0}}$,
and that $\GB(V)$ is of type~$\Res_{E_{0}/\QQ}\U_{n,E_{0}}$.

\begin{remark}
	\label{Hodge_theory_A}
	Let $A$ be an abelian surface over~$\CC$.
	Recall that $\HB^{2}(A)(1)$ has dimension~$6$.
	Let $H$ be the transcendental part of~$\HB^{2}(A)(1)$ and
	let $\rho$ denote the Picard number of~$A$, so that $\dim_{\QQ}(H) + \rho = 6$.
	If $A$ is simple, then the Albert classification of endomorphism algebras of
	abelian varieties states that $\End(A) \otimes \QQ$ can be one of the
	following:
	\begin{case}
		\item The field of rational numbers,~$\QQ$.
			In this case $\rho = 1$ and $\GB(H)$ is of type~$\SO_{5,\QQ}$.
		\item \label{hilbmod} A real quadratic extension~$F/\QQ$.
			In this case $\rho = 2$ and $\GB(H)$ is of type~$\SO_{4,\QQ}$.
			By exemple~3.2.2(a) of~\cite{norme},
			we see that $\map{\Nm_{F/\QQ}(\HH^{1}(A))}[hook]{\bigwedge^{2} \HH^{1}(A) \cong \HH^{2}(A)}$,
			where~$\Nm(\_)$ is the norm functor studied in~\cite{norme}.
			This norm map identifies $\Nm_{F/\QQ}(\HH^{1}(A))(1)$ with the transcendental part $H$.
			Observe that consequently the Hodge group~$\HdgB(\HH^{1}(A)) = \Res_{F/\QQ}\SL_{2,F}$
			is a $(2:1)$-cover of~$\GB(H)$.
		\item An indefinite quaternion algebra~$D/\QQ$.
			(This means that $D \otimes_{\QQ} \RR \cong M_{2}(\RR)$.)
			In this case $\rho = 3$ and $\GB(H)$ is of type~$\SO_{3,\QQ}$.
		\item \label{simplecm1} A \cm{} field~$E/\QQ$ of degree~$4$.
			In this case $\rho = 2$ and $\GB(H)$ is of type~$\Res_{E_{0}/\QQ}\U_{1,E_{0}}$.
	\end{case}
	(Note that the endomorphism algebra of~$A$ cannot be an imaginary quadratic field, by theorem~5 of~\cite{Sh63}.)
	If $A$ is isogenous to the product of two elliptic curves $Y_{1} \times Y_{2}$, then there are the following options:
	\begin{case}[resume]
		\item \label{nonsimpleFA} The elliptic curves are not isogenous, and neither of them is of \cm{}~type,
			in which case $\rho = 2$ and $\GB(H)$ is of type~$\SO_{4,\QQ}$.
			Indeed, $\HdgB(Y_{1})$ and~$\HdgB(Y_{2})$ are isomorphic to~$\SL_{2,\QQ}$.
			Note that $H = \HB^{2}(A)(1)^{\tra}$ is isomorphic to
			the exterior tensor product~$\big(\HB^{1}(Y_{1})
			\mathbin{\raisebox{.5pt}{\scalebox{.9}{$\boxtimes$}}} \HB^{1}(Y_{2})\big)(1)$,
			We find that $\GB(H)$ is the image of the canonical map $\map{\SL_{2,\QQ} \times \SL_{2,\QQ}}{\GL(H)}$.
			The kernel of this map is $\langle (-1,-1) \rangle$.
		\item The elliptic curves are not isogenous, one has endomorphism algebra~$\QQ$,
			and the other has~\cm{} by an imaginary quadratic extension~$E/\QQ$.
			In this case $\rho = 2$ and $\GB(H)$ is of type~$\U_{2,\QQ}$.
		\item \label{simplecm2} The elliptic curves are not isogenous,
			and $Y_{i}$ (for $i = 1,2$) has~\cm{} by an imaginary quadratic extension~$E_{i}/\QQ$.
			Observe that $E_{1} \not\cong E_{2}$, since $Y_{1}$ and~$Y_{2}$ are not isogenous.
			Let $E/\QQ$ be the compositum of $E_{1}$ and~$E_{2}$, which is a \cm~field of degree~$4$ over~$\QQ$.
			In this case $\rho = 2$ and $\GB(H)$ is of type~$\Res_{E_{0}/\QQ}\U_{1,E_{0}}$.
		\item The elliptic curves are isogenous and have trivial endomorphism algebra.
			In this case $\rho = 3$ and $\GB(H)$ is of type~$\SO_{3,\QQ}$.
		\item The elliptic curves are isogenous and have \cm~by an imaginary quadratic extension~$E/\QQ$.
			In this case $\rho = 4$ and $\GB(H)$ is of type~$\U_{1,\QQ}$.
	\end{case}
\end{remark}

\section{The main theorem: the Mumford--Tate conjecture for the product of an abelian surface and a K3~surface}
\label{Proof}

\paragraph{}
Let $K$ be a finitely generated subfield of~$\CC$.
Let $A$ be an abelian surface over~$K$, and
let $\HA$ denote the transcendental part of the motive~$\HH^{2}(A)(1)$.
(The Hodge structure~$H$ in \cref{Hodge_theory_A} is the Betti realisation~$\HB(\HA)$ of~$\HA$.)
Let $X$ be a K3~surface over~$K$, and
let $\HX$ denote the transcendental part of the motive~$\HH^{2}(X)(1)$.

Recall from \cref{FG} that we associated a field~$F$ and a group~$\scrG$ with every Hodge structure~$V$ of K3~type.
The important properties of $F$ and~$G$ are that
\begin{itemize}
	\item $\scrG^{\der}$ is an almost direct product of absolutely simple groups over $F$; and
	\item $\Res_{F/\QQ}\scrG$ is isogenous to $\GB(V)$.
\end{itemize}
Let $F_{A}$ and~$\scrG_{A}$ be the field and group associated with~$\HB(\HA)$ as in \cref{FG}.
Similarly, let $F_{X}$ and~$\scrG_{X}$ be the field and group associated with~$\HB(\HX)$.
Concretely, for $F_{A}$ this means that
\[
	F_{A} \cong
	\begin{cases}
		\End(A) \otimes \QQ & \text{in \cref{hilbmod} (so $F_{A}$ is \tr{} of degree~$2$)} \\
		E_{A,0} & \text{in \cref{simplecm1,simplecm2} (so $F_{A}$ is \tr{} of degree~$2$)} \\
		\QQ & \text{otherwise.}
	\end{cases}
\]
Let $E_{X}$ be the endomorphism algebra of $\HX$.
We summarise the notation for easy review during later parts of this section:
\begin{itemize}
	\item[$K$] finitely generated subfield of $\CC$
	\item[$A$] abelian surface over $K$
	\item[$\HA$] transcendental part of the motive $\HH^{2}(A)(1)$
	\item[$F_{A}$] field associated with the Hodge structure~$\HB(\HA)$, as in \cref{FG}
	\item[$\scrG_{A}$] group over $F_{A}$ such that $\Res_{F_{A}/\QQ}\scrG_{A}$ is isogenous to~$\GB(\HA)$, as in \cref{FG}
	\item[$X$] K3~surface over $K$
	\item[$\HX$] transcendental part of the motive $\HH^{2}(X)(1)$
	\item[$F_{X}$] field associated with the Hodge structure~$\HB(\HX)$, as in \cref{FG}
	\item[$\scrG_{X}$] group over $F_{X}$ such that $\Res_{F_{X}/\QQ}\scrG_{X}$ is isogenous to~$\GB(\HX)$, as in \cref{FG}
	\item[$E_{X}$] the endomorphism algebra of $\HX$
\end{itemize}

The proof of the main \namecref{MTCAxX} (\labelcref{MTCAxX}) will take the remainder of this article.
There are four main parts going into the proof, which are \cref{KSarg,diffEnd,diffLie,mtcnonsplitl}.
The \cref{reductions,reduce_der,nonsplitl,mtcnonsplitlcor} are small reductions and intermediate results.
Together \cref{KSarg,diffEnd,diffLie} deal with almost all combinations of abelian surfaces and K3~surfaces.
\Cref{mtcnonsplitl} is rather technical, and is the only place in the proof
where we use that $\HX$ really is a motive coming from a K3~surface.

\begin{lemma}
	\label{reductions}
	\begin{itemize}
		\item The Mumford--Tate conjecture for $\HH^{2}(A \times X)(1)$
			is equivalent to $\MTC(\HA \oplus \HX)$.
		\item The $\ell$-adic realisations of $\HA \oplus \HX$ form a compatible system of $\ell$-adic representations.
	\end{itemize}
	\begin{proof}
		The first claim follows from \cref{MTC_Tate_indep}.
		The $\Hl^{2}(A \times X)(1)$ form a compatible system of $\ell$-adic representations
		and we only remove Tate classes to obtain $\Hl(\HA \oplus \HX)$;
		hence the $\ell$-adic realisations of $\HA \oplus \HX$ also form a compatible system of $\ell$-adic representation.
	\end{proof}
\end{lemma}

\begin{lemma}
	\label{reduce_der}
	If for some prime~$\ell$, the natural morphism
	\[
		\map[\iota_{\ell}]{\Glc(\HA \oplus \HX)^{\der}}[hook]{\Glc(\HA)^{\der} \times \Glc(\HX)^{\der}}
	\]
	is an isomorphism, then the Mumford--Tate conjecture for $\HA \oplus \HX$ is true.
	\begin{proof}
		By \cref{MTCcentres}, we know that the Mumford--Tate conjecture for $\HA \oplus \HX$ is true on the centres
		of $\GB(\HA \oplus \HX) \otimes \QQl$ and $\Glc(\HA \oplus \HX)$.
		By Deligne's theorem on absolute Hodge cycles,
		we know that $\Glc(\HA \oplus \HX) \subset \GB(\HA \oplus \HX) \otimes \QQl$.
		Hence if $\map[\iota_{\ell}]{\Glc(\HA \oplus \HX)^{\der}}[hook]{\Glc(\HA)^{\der} \times \Glc(\HX)^{\der}}$
		is an isomorphism, then $\MTC_{\ell}(\HA \oplus \HX)$ is true,
		and by \cref{MTC_lindep}, so is $\MTC(\HA \oplus \HX)$.
	\end{proof}
\end{lemma}

\begin{lemma}
	\label{diffEnd}
	The Mumford--Tate conjecture for $\HA \oplus \HX$ is true if $F_{A} \not\cong F_{X}$.
	\begin{proof}
		By \cref{reduce_der} we are done if
		$\map[\iota_{\ell}]{\Glc(\HA \oplus \HX)^{\der}}[hook]{\Glc(\HA)^{\der} \times \Glc(\HX)^{\der}}$
		is an isomorphism for some prime~$\ell$.
		We proceed by assuming that for all~$\ell$, the morphism~$\iota_{\ell}$ is not an isomorphism.
		This will imply that $F_{A} \cong F_{X}$.

		By \cref{Galois_obs_ell},
		we see that $F_{A,\ell} = F_{A} \otimes \QQl$ and $F_{X,\ell} = F_{X} \otimes \QQl$ have an isomorphic factor,
		since we assume that $\iota_{\ell}$ is not an isomorphism.
		If $F_{A}$ is isomorphic to~$\QQ$, then $F_{X,\ell}$ has a factor~$\QQl$ for each~$\ell$,
		and we win by \cref{trans}.

		Next suppose that $F_{A} \not\cong \QQ$, in which case it is a real quadratic extension of~$\QQ$.
		In particular $F_{A}$ is Galois over~$\QQ$
		and $\scrG_{A}^{\der}$ is an absolutely simple group over~$F_{A}$ of Lie type~$A_{1}$.
		Using \cref{summand} we find, for each prime~$\ell$, semisimple Lie algebras
		$\mathfrak{s}_{A,\ell}$, $\mathfrak{t}_{\ell}$ and~$\mathfrak{s}_{X,\ell}$
		such that
		\begin{align*}
			\Lie(\scrG_{A}^{\der}) \cong
			\Lie(\Glc(\HA)^{\der}) &\cong \mathfrak{s}_{A,\ell} \oplus \mathfrak{t}_{\ell} \\
			\Lie(\scrG_{X}^{\der}) \cong
			\Lie(\Glc(\HX)^{\der}) &\cong
			\phantom{\mathfrak{s}_{A,\ell} \oplus \mathllap{\mathfrak{t}}} \mathfrak{t}_{\ell} \oplus \mathfrak{s}_{X,\ell} \\
			\Lie(\Glc(\HA \oplus \HX)^{\der}) &\cong \mathfrak{s}_{A,\ell} \oplus
			\mathfrak{t}_{\ell} \oplus \mathfrak{s}_{X,\ell}.
		\end{align*}
		The absolute ranks of these Lie algebras do not depend on~$\ell$,
		by \cref{MTCcentres} and remark~6.13 of~\cite{LP2} (or the letters of Serre to Ribet in~\cite{Serre}).

		If $\ell$ is a prime that is inert in $F_{A}$,
		then $\scrG_{A}^{\der} \otimes_{F_{A}} F_{A,\ell}$ is an absolutely simple group.
		Since $\mathfrak{t}_{\ell} \ne 0$, we conclude that $\mathfrak{s}_{A,\ell} = 0$.
		By the independence of the absolute ranks, $\mathfrak{s}_{A,\ell} = 0$ for all primes~$\ell$.
		Consequently, if $\ell$~is a prime that splits in~$F_{A}$,
		then $\mathfrak{t}_{\ell}$ has two simple factors that are absolutely simple Lie algebras over~$\QQl$ of Lie type $A_{1}$.

		If $\scrG_{X}^{\der}$ is an absolutely simple group over~$F_{X}$,
		then $F_{X,\ell}$ contains two copies of~$\QQl$, for each~$\ell$ that splits in~$F_{A}$.
		Recall that by \cref{Galois_obs_ell}, for all primes~$\ell$,
		we know that $F_{A,\ell}$ and $F_{X,\ell}$ have an isomorphic factor.
		In particular, for inert primes~$\ell$, $F_{A,\ell}$ is a factor of $F_{X,\ell}$.
		Hence $F_{A,\ell}$ is a factor of $F_{X,\ell}$ for all primes~$\ell$, and
		we are done, by \cref{trans}.

		If $\scrG_{X}^{\der}$ is not an absolutely simple group, then it is of type~$\SO_{4,E_{X}}$.
		In particular $\dim_{E_{X}}(\HX) = 4$ and $F_{X} \cong E_{X}$.
		It follows from \cref{Hodge_theory_K3} and the fact that $\dim_{\QQ}(\HX) \le 22$ that $[F_{X} : \QQ] \le 5$.
		By \cref{hack} we conclude that $F_{A} \cong F_{X}$.
	\end{proof}
\end{lemma}

\paragraph{}
\label{remaining}
From now on, we assume that $F_{A} \cong F_{X}$, which we will simply denote with~$F$.
We single out the following cases, and prove the Mumford--Tate
conjecture for $\HA \oplus \HX$ for all other cases in the next lemma.
\begin{case}
\item $\GB(\HA)$ and~$\GB(\HX)$ are both of type~$\SO_{5,\QQ}$;
\item $\GB(\HA)$ is of type $\SO_{3,\QQ}$, or~$\SO_{4,\QQ}$, or~$\U_{2,\QQ}$,
	and the type of~$\GB(\HX)$ is also one of these types;
\item \label{nasty} $F$ is a real quadratic extension of~$\QQ$, $A$~is an absolutely simple abelian surface with endomorphisms by~$F$
	(so $\scrG_{A} \cong \SL_{2,F}$), and
	\begin{case}
	\item \label{nastyeasy} $\scrG_{X}$ is of type~$\SO_{3,F}$ or~$\U_{2,F}$; or
	\item \label{nastyhard} $\scrG_{X}$ is non-simple of type~$\SO_{4,F}$ as in \cref{SO4nonsimple} of~\cref{SO4}.
	\end{case}
\end{case}
We point out that in the first two cases $\dim(\HX) \le 5$, which can be deduced from \cref{Hodge_theory_K3}.

\begin{lemma}
	\label{diffLie}
	If we are not in one of cases listed in \cref{remaining}, then
	the Mumford--Tate conjecture for $\HA \oplus \HX$ is true.
	\begin{proof}
		By \cref{reduce_der} we are done if
		$\map[\iota_{\ell}]{\Glc(\HA \oplus \HX)^{\der}}[hook]{\Glc(\HA)^{\der} \times \Glc(\HX)^{\der}}$
		is an isomorphism for some prime~$\ell$.

		The crucial ingredient in this lemma is \cref{Lie_obs}.
		Recall that $\CC \cong \widebar{\QQl}\mkern2mu$, as fields.
		If the Dynkin diagram of $\Lie(\Glc(\HA)^{\der})_{\CC}$ has no components in common
		with the Dynkin diagram of $\Lie(\Glc(\HX)^{\der})_{\CC}$,
		by \cref{Lie_obs}, we see that $\iota_{\ell}$~is an isomorphism, and we win.
		Recall that $\MTC(\HA)$ and~$\MTC(\HX)$ are known.
		Thus $\iota_{\ell}$~is an isomorphism when
		the Dynkin diagram of $\Lie(\GB(\HA)^{\der})_{\CC}$ has no components in common
		with the Dynkin diagram of $\Lie(\GB(\HX)^{\der})_{\CC}$.
		By inspection of \cref{Hodge_theory_K3} and \cref{Hodge_theory_A},
		we see that this holds, except for the cases listed in~\cref{remaining}.
	\end{proof}
\end{lemma}

\begin{lemma}
	\label{KSarg}
	The Mumford--Tate conjecture for $\HA \oplus \HX$ is true if $\dim(\HX) \le 5$.
	In particular, the Mumford--Tate conjecture is true for the first two cases listed in \cref{remaining}.
	\begin{proof}
		Let $B$ be the Kuga--Satake variety associated with~$\HB(\HX)$.
		This is a complex abelian variety of dimension $2^{\dim(\HX) - 2}$.
		Up to a finitely generated extension of~$K$, we may assume that $B$ is defined over~$K$.
		(In fact, $B$~\emph{is} defined over~$K$, by work of Rizov, \cite{Rizov}.)
		By \cref{fgext}, we may and do allow ourselves a finite extension of~$K$,
		to assure that $B$ is isogenous to a product of absolutely simple abelian varieties over~$K$.
		By proposition~6.3.3 of~\cite{vG00}, we know that $\HB(\HX)$ is a sub-$\QQ$-Hodge structure of~$\End(\HB^{1}(B))$.
		Since $\HX$ is an abelian motive, we deduce that $\HX$ is a submotive of~$\End(\HH^{1}(B))$,
		by Andr\'{e}'s ``Hodge = motivated'' theorem (see th\'{e}ori\`{e}me~0.6.2 of~\cite{An95}).
		Consequently, $\MTC(A \times B)$ implies $\MTC(\HA \oplus \HX)$.

		Recall that the even Clifford algebra~$C^{+}(\HX) = C^{+}(\HB(\HX))$ acts on~$B$.
		Theorem~7.7 of~\cite{vG00} gives a description of~$C^{+}(\HX)$;
		thus describing a subalgebra of~$\End^{0}(B)$.
		\begin{itemize}
			\item If $\dim(\HX) = 3$, then $\dim(B) = 2$ and
				$C^{+}(\HX)$ is a quaternion algebra over~$\QQ$.
			\item If $\dim(\HX) = 4$, then $\dim(B) = 4$ and
				$C^{+}(\HX)$ is either a product $D \times D$,
				where $D$ is a quaternion algebra over~$\QQ$;
				or $C^{+}(\HX)$ is a quaternion algebra over a
				totally real quadratic extension of~$\QQ$.
			\item If $\dim(\HX) = 5$, then $\dim(B) = 8$ and
				$C^{+}(\HX)$ is a matrix algebra~$M_{2}(D)$,
				where $D$ is a quaternion algebra over~$\QQ$.
		\end{itemize}
		We claim that $A \times B$ satisfies the conditions of \cref{MTCgLt}.
		First of all, observe that $A$ satisfies those conditions,
		which can easily be seen by reviewing \cref{Hodge_theory_A}.
		We are done if we check that $B$ satisfies the conditions as well.
		\begin{itemize}
			\item If $\dim(\HX) = 3$, then $B$ is either a simple
				abelian surface, or isogenous to the square of
				an elliptic curve. In both cases, $B$ satisfies
				the conditions of \cref{MTCgLt}.
			\item If $\dim(\HX) = 4$, and $C^{+}(\HX)$ is $D \times
				D$ for some quaternion algebra~$D$ over~$\QQ$,
				then $B$ splits (up to isogeny) as $B_{1} \times
				B_{2}$. In particular $\dim(B_{i}) = 2$, since
				$D$ cannot be the endomorphism algebra of an
				elliptic curve. Hence both $B_{i}$ satisfy the
				conditions of \cref{MTCgLt}.

				On the other hand, if $\dim(\HX) = 4$ and
				$C^{+}(\HX)$ is a quaternion algebra over a
				totally real quadratic extension of~$\QQ$, then
				there are two options.
				\begin{itemize}
					\item If $B$ is not absolutely simple,
						then all simple factors have
						dimension $\le 2$; since
						$\End^{0}(B)$ is
						non-commutative.
						Indeed, the product of an
						elliptic curve and a simple
						abelian threefold has
						commutative endomorphism ring
						(see, \textit{e.g.,} section~2
						of~\cite{MZ99}).
					\item If $B$ is absolutely simple, then it has
						relative dimension~$1$. This
						abelian fourfold must be of
						type~\textsc{ii}$(2)$, since
						type~\textsc{iii}$(2)$ does not
						occur (see proposition~15
						of~\cite{Sh63}, or table~1
						of~\cite{MZ95} which also
						proves~$\MTC(B)$).
				\end{itemize}
				In both of these cases, $B$ satisfies the
				conditions of \cref{MTCgLt}.
			\item If $\dim(\HX) = 5$, then $B$ is the square of an
				abelian fourfold~$C$ with endomorphism algebra
				containing a quaternion algebra over~$\QQ$.
				\begin{itemize}
					\item If $C$ is not absolutely simple,
						then all simple factors have
						dimension $\le 2$; since
						$\End^{0}(C)$ is
						non-commutative.
					\item If $C$ is simple, then we claim
						that $C$ must be of type~\textsc{ii}.
						Indeed, since $\HB(\HX)$ is a
						sub-$\QQ$-Hodge structure
						of~$\End(\HB^{1}(B))$, the
						Mumford--Tate group of~$B$ must
						surject onto~$\GB(\HX)$. In this
						case, $\dim(\HX) = 5$, hence
						$\GB(\HX)$ is of type~$\SO_{5,\QQ}$,
						with Lie type~$B_{2}$.
						But \S6.1~of~\cite{MZ95}
						shows that if $C$ is of
						type~\textsc{iii}, then
						$\GB(C)$ has Lie type~$D_{2} \cong A_{1} \oplus A_{1}$.
						This proves our claim.
						Since $\End^{0}(C)$ is a
						quaternion algebra and $C$ is
						an abelian fourfold,
						table~1 of~\cite{MZ95} shows that
						$\MTC(C)$ is true
						and $D_{4}$ does not occur in
						the Lie type of~$\GB(C)$. 
				\end{itemize}
		\end{itemize}
		We conclude that $\MTC(A \times B)$ is true,
		and therefore $\MTC(\HA \oplus \HX)$ is true as well.
	\end{proof}
\end{lemma}

The only cases left are those listed in~\cref{nasty} of~\cref{remaining}.
Therefore, we may and do assume that $F$ is a real quadratic field extension of~$\QQ$;
and that $A$ is an absolutely simple abelian surface with endomorphisms by~$F$ (\textit{i.e.,} \cref{hilbmod}).
In particular $\scrG_{A} = \SL_{2,F}$.

\begin{lemma}
	\label{nonsplitl}
	If $X$ falls in one of the subcases listed in \cref{nasty},
	then there exists a place~$\lambda$ of~$F$
	such that $\scrG_{X}^{\der} \otimes_{F} F_{\lambda}$ does not contain a split factor.
	\begin{proof}
		In \cref{nastyeasy}, $\scrG_{X}^{\der}$ is of Lie type~$A_{1}$.
		In \cref{nastyhard}, $\scrG_{X} \sim N \times N^{\textnormal{op}}$,
		where $N$ is a form of~$\SL_{2,F}$, as explained in \cref{SO4}.
		By theorem~26.9 of~\cite{BoI}, there is an equivalence between forms of~$\SL_{2}$ over a field,
		and quaternion algebras over the same field.
		We find a quaternion algebra~$D$ over~$F$ corresponding to~$\scrG_{X}^{\der}$, respectively~$N$,
		in \cref{nastyeasy}, respectively \cref{nastyhard}.
		In particular $\scrG_{X}^{\der}$ contains a split factor if and only if the quaternion algebra is split.

		Let $\Set{\sigma,\tau}$ be the set of embeddings~$\Hom(F, \RR)$.
		Since $F$ acts on~$\HB(\HX)$, we see that
		$F \otimes_{\QQ} \RR \cong \RR^{(\sigma)} \oplus \RR^{(\tau)}$
		acts on
		\[
			\HB(\HX) \otimes_{\QQ} \RR \cong W^{(\sigma)} \oplus W^{(\tau)}.
		\]
		Here $W^{(\sigma)}$ and~$W^{(\tau)}$ are $\RR$-Hodge structure of dimension~$\dim_{F}(\HX)$.
		Observe that the polarisation form is definite on one of the terms, while it is non-definite on the other.
		Without loss of generality we may assume that the polarisation form is definite on $W^{(\sigma)}$,
		and non-definite on $W^{(\tau)}$.

		Thus, the group~$\GB(\HX) \otimes_{\QQ} \RR$ is the product of a compact group and a non-compact group;
		and therefore, $\Res_{F/\QQ}\scrG_{X} \otimes_{\QQ} \RR$ is the product of a compact group and a non-compact group.
		Indeed $\scrG_{X} \otimes_{F} \RR^{(\sigma)}$ is compact,
		while $\scrG_{X} \otimes_{F} \RR^{(\tau)}$ is non-compact.
		By the first paragraph of the proof, this means that $D \otimes_{F} \RR^{(\sigma)}$ is non-split,
		while $D \otimes_{F} \RR^{(\tau)}$ is split.

		Since the Brauer invariants of~$D$ at the infinite places do not add up to~$0$,
		there must be a finite place~$\lambda$ of~$F$ such that $D_{\lambda}$ is non-split.
		Therefore $\scrG_{X}^{\der} \otimes_{F} F_{\lambda}$ does not contain a split factor.
	\end{proof}
\end{lemma}

\begin{lemma}
	\label{mtcnonsplitl}
	Assume that $K$ is a number field.
	If $X$ falls in one of the subcases listed in~\cref{nasty},
	then there is a prime number $\ell$ for which the natural map
	\[
		\map[\iota_{\ell}]{\Glc(\HA \oplus \HX)^{\der}}[hook]{\Glc(\HA)^{\der} \times \Glc(\HX)^{\der}}
	\]
	is an isomorphism.
	\begin{proof}
		The absolute rank of $\Glc(\HA \oplus \HX)^{\der}$ does not depend on~$\ell$,
		by \cref{MTCcentres,reductions} and
		remark~6.13 of~\cite{LP2} (or the letters of Serre to Ribet in~\cite{Serre}).
		Let $\ell$ be a prime that is inert in~$F$.
		Observe that all simple factors of $\Lie(\Glc(\HA)^{\der} \times \Glc(\HX)^{\der})$
		are $\QQl$-Lie algebras with even absolute rank (since $[F : \QQ] = 2$).
		By \cref{summand}, the Lie algebra of $\Glc(\HA \oplus \HX)^{\der}$
		is a summand of $\Lie(\Glc(\HA)^{\der} \times \Glc(\HX)^{\der})$,
		and therefore the absolute rank of $\Glc(\HA \oplus \HX)^{\der}$ must be even.

		Let $\lambda$ be one of the places of~$F$ found in \cref{nonsplitl},
		and let $\ell$ be the place of~$\QQ$ lying below~$\lambda$.
		Since $\Lie(\Glc(\HA \oplus \HX)^{\der})$ must surject to~$\Lie(\Glc(\HA))$
		(which is split, and has absolute rank~$2$),
		and $\Lie(\Glc(\HA \oplus \HX)^{\der})$ must also surject onto $\Lie(\Glc(\HX)^{\der})$,
		which has no split factor, by \cref{nonsplitl},
		we conclude that the absolute rank of $\Lie(\Glc(\HA \oplus \HX)^{\der})$ must be at least~$3$.
		By the previous paragraph, we find that the absolute rank must be at least $4$.

		If $\dim_{E_{X}}(\HX) \ne 4$ (\cref{nastyeasy})
		then $\scrG_{X}^{\der}$ is a group of Lie type~$A_{1}$,
		and therefore the product $\Glc(\HA)^{\der} \times \Glc(\HX)^{\der}$ has absolute rank~$4$.
		Hence $\Glc(\HA \oplus \HX)^{\der}$ must have absolute rank~$4$,
		which means that $\iota_{\ell}$ is an isomorphism, by \cref{summand,isomLie}.

		If $\dim_{E_{X}}(\HX) = 4$ (\cref{nastyhard}),
		then $\scrG_{X}$ is a group of Lie type $D_{2} = A_{1} \oplus A_{1}$.
		(Note that in this final case $\GB(\HA)$ and~$\GB(\HX)$ are semisimple,
		and therefore we may drop all the superscripts $(\_)^{\der}$ from the notation.)
		Since in this case $\Glc(\HA) \times \Glc(\HX)$ has absolute rank~$6$,
		and the absolute rank of $\Glc(\HA \oplus \HX)$
		is~$\ge 4$, it must be $4$~or~$6$ (since it is even).

		Suppose $\Glc(\HA \oplus \HX)$ has absolute rank~$4$.
		We apply \cref{summand} to the current situation, and find Lie algebras
		$\mathfrak{t}$ and~$\mathfrak{s}_{2}$ over~$\QQl$
		such that $\Lie(\Glc(\HA)) \cong \mathfrak{t}$
		and~$\Lie(\Glc(\HA \oplus \HX)) \cong \Lie(\Glc(\HX)) \cong \mathfrak{t} \oplus \mathfrak{s}_{2}$.
		In particular, $\Lie(\Glc(\HX))$ which is
		isomorphic to~$\Lie(\scrG_{X}) \otimes \QQl$ has a split simple factor.
		By \cref{nonsplitl} this means that $\ell$~splits in~$F$ as $\lambda \cdot \lambda'$.
		Observe that $F_{\lambda} \cong \QQl \cong F_{\lambda'}$.

		Note that in this case $\GB(\HX) \cong \Res_{F_{X}/\QQ} \scrG_{X}$,
		and since $\MTC(\HX)$ is known
		we find $\Glc(\HX) \cong \scrG_{X,\lambda} \times \scrG_{X,\lambda'}$ and
		a decomposition $\Hl(\HX) \cong \HH_{\lambda}(\HX) \oplus \HH_{\lambda'}(\HX)$.
		The group~$\Gal(\widebar{K}/K)$ acts on~$\HH_{\lambda}(\HX)$ via~$\scrG_{X,\lambda}$,
		and on~$\HH_{\lambda'}(\HX)$ via~$\scrG_{X,\lambda'}$.

		To summarise, our situation is now as follows.
		The prime number~$\ell$ splits in~$F$ as $\lambda \cdot \lambda'$.
		The group~$\scrG_{A}$ is isomorphic to~$\SL_{2,F}$, and is split and simply connected,
		The group~$\scrG_{X,\lambda'}$ is split, of type~$\SO_{4,\QQl}$, with Lie algebra $\mathfrak{t}$.
		The group~$\scrG_{X,\lambda}$ is non-split, of type~$\SO_{4,\QQl}$, with Lie algebra $\mathfrak{s}_{2}$.
		Recall the natural diagram:
		\[
			\begin{tikzpicture}[commutative diagrams/every diagram]
				\node (GlHAHX) at (0,3)
				{ $\Glc(\HA \oplus \HX)$ };
				\node (GlHAxHX) at (0,1.5)
				{ $\Glc(\HA) \times \Glc(\HX)$ };
				\node (GlHA) at (-3,0)
				{ $(\SL_{2,\QQl} \times \SL_{2,\QQl})/\langle (-1,-1) \rangle \cong \Glc(\HA)$ };
				\node (GlHX) at (3,0)
				{ $\Glc(\HX) \cong \scrG_{X,\lambda'} \times \scrG_{X,\lambda}$ };

				\path[commutative diagrams/.cd, every arrow, every label]
				(GlHAHX) edge node[commutative diagrams/hook] {$\iota_{\ell}$} (GlHAxHX)
				(GlHAxHX) edge node[commutative diagrams/two heads] {} (GlHA)
				(GlHAxHX) edge node[commutative diagrams/two heads] {} (GlHX);
			\end{tikzpicture}
		\]

		We are now set for the attack.
		We claim that the Galois representations $\Hl(\HA)$ and~$\HH_{\lambda'}(\HX)$ are isomorphic.
		Indeed, from the previous paragraph we conclude that $\Glc(\HA \oplus \HX) \cong \Gamma \times \scrG_{X,\lambda}$,
		where $\Gamma$~is a subgroup of $\Glc(\HA) \times \scrG_{X,\lambda'}$ with surjective projections.
		Thus $\Hl(\HA)$ and~$\HH_{\lambda'}(\HX)$ are both orthogonal representations of~$\Gal(\widebar{K}/K)$,
		and the action of Galois factors via~$\Gamma(\QQl)$.

		The Lie algebra of~$\Gamma$ is isomorphic to $\mathfrak{t}$,
		and $\Lie(\Gamma)$~is the graph of an isomorphism $\map{\Lie(\Glc(\HA)}{\Lie(\scrG_{X,\lambda'})}$.
		Since $\Glc(\HA)$ and~$\scrG_{X,\lambda'}$ have
		$(2:1)$-covers by~$\Res_{F/\QQ}\SL_{2,F} \cong \Hdgl(A)$ with kernels~$\Set{ \pm 1 }$,
		and $\Gamma$~is a subgroup of~$\Glc(\HA) \times \scrG_{X,\lambda'}$,
		we find that $\Gamma$ also has a $(2:1)$-cover by~$\Res_{F/\QQ}\SL_{2,F}$.
		Hence $\Gamma$ is the graph of an isomorphism~$\map{\Glc(\HA)}{\scrG_{X,\lambda'}}$.
		Because $\Hl(\HA)$ and~$\HH_{\lambda'}(\HX)$ are $4$-dimensional faithful orthogonal representations of $\Gamma$,
		they must be isomorphic; for up to isomorphism, there is a unique such representation.

		As a consequence, for any place~$v$ of~$K$,
		the characteristic polynomial of~$\Frob_{v}$ acting on~$\Hl(\HA)$
		coincides with its characteristic polynomial when acting on~$\HH_{\lambda'}(\HX)$.
		We conclude that $\charpol_{F_{\lambda'}}(\Frob_{v} | \HH_{\lambda'}(\HX))$ has coefficients in~$\QQ$.
		But then the same is true for $\charpol_{F_{\lambda}}(\Frob_{v} | \HH_{\lambda}(\HX))$
		since their product is $\charpol_{\QQl}(\Frob_{v} | \Hl(\HX))$, which has coefficients in~$\QQ$.

		Since we assumed that $K$ is a number field, we may apply the following results:
		\begin{itemize}
			\item Theorem~1 (item~1) of \cite{FC14}, which tells us that
				(up to a finite extension of~$K$, which does not matter, by \cref{fgext})
				there exists a set~$\mathscr{V}$ of places of~$K$ with density~$1$
				such that $X$ has good reduction at places $v \in \mathscr{V}$,
				and the Picard number of the reduction~$X_{v}$
				is the same as that of~$X$ (which, in our case is $22-8=14$).
			\item Proposition~3.2 of~\cite{YuYui}, which says that if $X$ has good and ordinary reduction at~$v$,
			      then the characteristic polynomial $\charpol_{\QQl}(\Frob_{v} | \Hl^{2}(X_{v})^{\tra})$ is
			      a power of an irreducible polynomial with coefficients in~$\QQ$.
		\end{itemize}
		We find that
		$\charpol_{F_{\lambda'}}(\Frob_{v} | \HH_{\lambda'}(\HX)) =
		\charpol_{F_{\lambda}}(\Frob_{v} | \HH_{\lambda}(\HX))$,
		for all places $v \in \mathscr{V}$.
		Since $\Gal(\widebar{K}/K)$ is compact, we may apply the argument given on the first pages of~\cite{Cours},
		and find that $\HH_{\lambda'}(\HX) \cong \HH_{\lambda}(\HX)$ as Galois representations.
		This contradicts the fact that $\scrG_{X,\lambda'}$ is split, while $\scrG_{X,\lambda}$ is not.
		We conclude that the rank must be~$6$, which implies, by \cref{summand,isomLie},
		that~$\iota_{\ell}$ is an isomorphism.
	\end{proof}
\end{lemma}

\begin{corollary}
	\label{mtcnonsplitlcor}
	If $X$ falls in one of the subcases listed in~\cref{nasty},
	then the Mumford--Tate conjecture is true for $\HA \oplus \HX$.
	\begin{proof}
		By \cref{reduce_der} we are done if
		$\map[\iota_{\ell}]{\Glc(\HA \oplus \HX)^{\der}}[hook]{\Glc(\HA)^{\der} \times \Glc(\HX)^{\der}}$
		is an isomorphism for some prime~$\ell$.
		This result follows from \cref{specialise,mtcnonsplitl}.
	\end{proof}
\end{corollary}

\paragraph{Proof of \cref{MTCAxX}}
We know have all tools in place to prove the main theorem.
% \begin{proof}[of \cref{MTCAxX}]
	By \cref{reductions} we reduce to the Mumford--Tate conjecture for $\HA \oplus \HX$.
	The theorem follows from \cref{KSarg,diffEnd,diffLie,mtcnonsplitlcor}. \quelle{$\square$}
% \end{proof}

% \begin{remark}
% 	By \cref{MTCAxX} we know that $\HB(\HA \otimes \HX)$ contains a non-trivial Hodge class
% 	if and only if $\Hl(\HA \otimes \HX)$ contains a non-trivial Tate class.
% 	In particular, if $\HB(\HA \otimes \HX)$ does not contain a non-trivial Hodge class, then
% 	the Hodge conjecture is true for~$\HH^{4}(A \times X)(2)$; and the Tate conjecture follows from \cref{MTCAxX}.

% 	If $\HB(\HA \otimes \HX)$ contains a non-trivial Hodge class $\gamma$,
% 	then~$\HB(\HA)$ and~$\HB(\HX)$ are isomorphic (since they are self-dual and irreducible).
% 	Because $\HA$ and~$\HX$ are abelian motives, this implies that $\HA$ and~$\HX$ are isomorphic motives.
% 	By transport of structure $\HA$ and~$\HX$ must be isometric for any choice of polaristation form on them.
% 	If furthermore $\HA$ and~$\HX$ have an endomorphism algebra that is a \cm~field,
% 	then recent work of Buskin (see theorem~1.1 of~\cite{Buskin})
% 	shows that there is an algebraic correspondence~$\Gamma$ between $X$ and~$A$
% 	whose cycle class $\operatorname{cl}(\Gamma) \in \HB^{4}(A \times X)(2)$ projects onto $\gamma \in \HB(\HA \otimes \HX)$.
% 	In particular, the Hodge conjecture is true for~$\HH^{4}(A \times X)(2)$; and the Tate conjecture follows from \cref{MTCAxX}.
% 	However, if the endomorphism algebra of~$\HA$ and~$\HX$ is a \tr~field, then there is nothing we can say.
% \end{remark}

\printbibliography

\end{document}